\documentclass[11pt]{article}
\usepackage{latexsym,amsmath}
\usepackage{fullpage}
\usepackage{dsfont}
\long\def\ignore#1{}

\newcommand\comm[1]{\typeout{Used \string\comm...}\quad{\bf [[#1]]}\quad}

\newtheorem{theorem}{Theorem}

\newtheorem{lemma}[theorem]{Lemma}

\def\qed{\hfill$\Box$\medskip}
\def\N{{\mathds N}}
\def\R{{\mathds R}}

\def\eps{\varepsilon}

\title{On a question of Bourgain about geometric incidences\thanks{The
research was partially supported by OTKA and NSERC grants. The first
author was supported by a Sloan Research Fellowship.}}

\author{
J\'ozsef Solymosi\thanks{Department of Mathematics, University of
British Columbia, Vancouver, BC, Canada,
\texttt{solymosi@math.ubc.ca}} \and Csaba D.
T\'oth\thanks{Department of Mathematics, Room 2-336, MIT, Cambridge,
MA~02139, USA, \texttt{toth@math.mit.edu}}}

\date{}

\begin{document}
\maketitle

\begin{abstract}
  Given a set of $s$ points and a set of $n^2$ lines in
  three-dimensional Euclidean space such that each line is incident to
  $n$ points but no $n$ lines are coplanar, then we have
  $s=\Omega(n^{11/4})$. This is the first nontrivial answer to a
  question recently posed by Jean Bourgain.
\end{abstract}

\section{Introduction\label{sec:intro}}

A seminal result in geometric incidences is the Szemer\'edi-Trotter
bound~\cite{st-epdg-83}, which says that the number of incidences
between $n$ points and $m$ lines in the Euclidean plane is at most
$O(n^{2/3}m^{2/3}+n+m)$, and it is the best possible. Typical problems
in geometric incidences consider two families of geometric objects of
size $n$ and $m$, respectively, and asks for the maximum number of
incident pairs.
In three and higher dimensions, non-trivial bounds often require
restrictions on the geometric objects, otherwise every object of one
family might be incident to all objects of the other family, or the
number of incidences might be maximized when all objects lie in a
lower dimensional affine subspace.

In the 1990-s, Tom Wolff~\cite{W1,W2} observed that tight bounds on
the number of incidences can be used efficiently to attack problems
related to the Kakeya conjecture, one of the central conjectures in
harmonic analysis. For more details about the Kakeya conjecture, refer
to \cite{b-sttsp-99} or \cite{W2}.
Bennett, Carbery, and Tao~\cite{BCT} established a concrete connection
between multilinear Kakeya estimates and bonds on number of incidences
between points and lines in three dimensions.

Recently, Jean Bourgain~\cite{cl-ppwrt-04} asked what is the minimum
cardinality $s$ of a point set $S$ in three-dimensional Euclidean
space, if we are given $n^2$ lines, each of which1 is incident to $n$
points of $S$ but no $n$ lines are coplanar.  He conjectured that
$s=\Omega(n^{3-\eps})$ for every constant $\eps>0$.  If the $n^2$
lines are disjoint, then $s=n^3$ is obvious.  The integer grid
$\{(a,b,c)\in \R^3 :1\leq a,b,c\leq n\}$ and the set of $3n^2$
axis-aligned lines also gives $s=n^3$ points.
The Szemer\'edi-Trotter theorem, applied to $n^2$ lines and $s$
points with $n^3$ point-line incidences, gives a lower bound of
$s=\Omega(n^{5/2})$. This bound, however, does not use the condition
that no $n$ lines are coplanar. We give the first nontrivial answer
to Bourgain's question.

\begin{theorem}\label{thm:main}
  Given a set of $s$ points and a set of $n^2$ lines in $\R^3$ such
  that every line is incident to at least $n$ points but no $n$ lines
  are coplanar, we have $s=\Omega(n^{11/4})=\Omega(n^{3-\frac{1}{4}})$.
\end{theorem}

\paragraph{Related previous work.}
The Szemer\'edi-Trotter bound on the number of point-line incidences
is tight in the plane: $n$ points and $m$ lines in the plane may have
$\Omega(n^{2/3}m^{2/3}+n+m)$ incidences. A better upper bound in
three-dimensions is possible only under certain restrictions, which
guarantee that the given point-line configuration is ``far'' from
being planar. Bourgain's condition that no $n$ of $n^2$ given lines
are coplanar is one such restriction.  Another previously considered
condition requires every point to be incident to at least three
non-coplanar lines, such a point is called {\em joint}.
Sharir~\cite{s-ojals-94} conjectured that the number of joints for $m$
lines in $\R^3$ is $O(m^{3/2})$, which is attained by the axis-aligned
lines of the $\sqrt{m}\times \sqrt{m}\times \sqrt{m}$ integer lattice
section. The best current bound, $O(m^{1.6232})$ by Feldman and
Sharir~\cite{fs-ibjal-05}, is far from being tight. Sharir and
Welzl~\cite{sw-plis-04} gave an $O(m^{5/3})$ upper bound on the number
of joint-line incidences in three-space.

Edelsbrunner {\em et al.}~\cite{egs-cmcap-90} were the first to study
the number of point-plane incidences in three dimensions. They
obtained an $O(n^{3/5}m^{4/5} +n +m\log n)$ bound for $n$ points and
$m$ planes assuming that no three points are collinear. Note that $n$
points and $m$ planes in three-space may have $nm$ incidences if the
points are collinear and all $m$ planes contain this line. Every
non-trivial bound on point-plane incidences must, therefore, impose
some reasonable restriction.
Agarwal and Aronov~\cite{aa-cfi-92} proved an upper bound of $O(n
m^{2/3} + n^2)$ for point-plane incidences in $\R^3$ assuming that
each plane is spanned by the point set (that is, each plane contains
three affine independent points). Their bound matches the lower bound
of Edelsbrunner and Haussler~\cite{eh-cctda-86}.
Bra\ss\ and Knauer~\cite{bk-ocphi-03} gave an $O(n^{3/4} m^{3/4}
\log (nm) + (n+m) \log (n+m))$ bound assuming that the incidence
graph does not contain a $K_{r,r}$ for some fixed $r\in \N$.
Elekes and T\'oth~\cite{et-inh-05} obtained a tight bound of
$O(n^{3/4}m^{3/4}+n\sqrt{m}+n)$ for the incidences between $n$ points
and $m$ {\em saturated} planes (where a plane is called saturated if
at most a constant fraction of the points lying in the plane are
collinear).  It is attained by a system where all points lie in two
parallel planes.
Solymosi and T\'oth~\cite{st-ddhse-06} gave an $O(n^{3/4}m^{3/4})$
bound for homogeneous point sets, which covers the example of the
integer lattice section.
Solymosi and Vu \cite{sv-ddhs-03} shoed that if $S$ is a homogeneous
set of $n$ points in three-space and $k\geq 2$, then the number of
$k$-rich lines is at most $O(n^2/k^4)$. Note that this result implies
Bourgain's conjecture for homogeneous point sets.

\paragraph{Proof techniques and Organization.}
Essentially two different methods have been developed for proving
geometric incidence bounds: One is the {\em crossing number}
technique based on work by Sz\'ekely~\cite{s-cnhep-97}; the other is
the {\em $\eps$-cutting} technique, which is a divide-and-conquer
strategy introduced by Clarkson and Shor~\cite{cs-arscg-89}, and
some tight bounds were obtained by Chazelle and
Friedman~\cite{cf-dvrsu-90} (see also, \cite{c-c-05,matousek}). We
deploy both techniques. Refer to a survey by Pach and
Sharir~\cite{ps-gi-03} for the rich history and widespread
applications of these techniques.

In Section~\ref{sec:reg}, we use ideas of Sharir and
Welzl~\cite{sw-plis-04} to represent lines meeting a given line
and doubly ruled surfaces by points and algebraic curves in the
plane; and then we apply the crossing technique in the plane.  In
Section~\ref{sec:proof}, we apply the cutting technique to reduce the
problem to the case that every point is incident to $\Omega(1)$ lines
on average, and complete the proof of Theorem~\ref{thm:main} with an
extremal graph theoretical bound on the number of incidences of lines
and doubly ruled surfaces.

\section{Lines and reguli in three-space\label{sec:reg}}

A {\em regulus} is a doubly ruled quadratic surface in
three-space~\cite{h-52}. Every regulus contains two families of lines,
which are called {\em ruling}s: Each ruling consists of pairwise skew
lines, and every line of one ruling intersects all the lines of the
other ruling of a regulus. We say that every line of each ruling is
{\em incident} to the regulus. Any three pairwise skew lines are
contained in a ruling of a unique regulus. If a line has three common
points with a regulus, then it must be incident to that regulus
(furthermore, it is part of one ruling of that regulus).

Since no three lines are incident to two distinct rulings, we can
apply the K\H{o}v\'ari-S\'os-Tur\'an bound from extremal graph
theory~\cite{kst-54}: A bipartite graph with $m$ and $r$ vertices in
its two vertex classes and having no subgraph isomorphic to $K_{3,2}$
has at most $O(mr^{2/3}+r)$ edges. It follows that the number of
line-regulus incidences between $m$ lines and $r$ reguli in
three-space is bounded by $O(mr^{2/3}+r)$.

We extend two lemmas of Sharir and Welzl~\cite{sw-plis-04} on the
number of point-regulus incidences. Both concern the number of reguli
incident to a given line $\ell$, and spanned by three lines of a given
set of lines $M$.  The first lemma gives an upper bound on the number
line-reguli incidences; the second lemma gives a lower bound on the
the number of reguli under the condition that at most $|M|/5$ lines of
$M$ may be coplanar or concurrent.

\begin{lemma}\label{lem:k/m}
Assume that a line $\ell$ meets every element of a set $M$ of $m$
lines in three-space and we are given a set $R$ of $r$ reguli. The
number of incidences between $R$ and $M$ is bounded by
$O(m^{3/5}r^{4/5}+m+r).$
\end{lemma}

\begin{proof}
  If a regulus $\varrho \in R$ is incident to two intersecting lines
  of $M$, it cannot be incident to line $\ell$.  Hence such a
  $\varrho$ intersects $\ell$ in at most two points, and so it
  contains at most four lines of $M$. The reguli in $R$ that are
  incident to up to four lines in $M$ are responsible for at most $4r$
  incidences. Let $Q$ denote the set of reguli in $R$ incident to at
  least three pairwise skew lines of $M$. We can represent the lines
  that meet $\ell$ by points in three-space: For instance, one
  dimension can be the intersection with $\ell$, and two additional
  dimensions can be the coordinates of their intersection point with a
  plane parallel to $\ell$. The families of lines incident to reguli
  spanned by $M$ correspond to bounded degree algebraic curves in
  three-space~\cite{sw-plis-04}.
  
  Project these points and curves to a generic plane.  We obtain a set
  $P$ of $m$ points and a set $C$ of at most $r$ bounded degree
  algebraic curves such that any three points of $P$ are incident to
  at most one curve of $C$. By result of Pach and
  Sharir~\cite{ps-onibp-98}, the number of point-curve incidences is
  bounded by $O(m^{3/5}r^{4/5}+m+r)$. Together with $O(r)$ incidences
  of the reguli in $R\setminus Q$, there are $O(m^{3/5}r^{4/5}+m+r)$
  line-regulus incidences between $M$ and $R$.  \qed
\end{proof}

\noindent{\bf Remark:}
Lemma~\ref{lem:k/m} is not sharp. With a little work, one could show
that the curves in the plane are pseudo-parabolas, for which a better
incidence bound is available~\cite{AAT}. Using this bound, one could
show that there are at most $O(m^{6/11-\epsilon}r^{9/11}+m+r)$
line-regulus incidences between $M$ and $R$. This bound, however, is
not a bottle neck in our estimates, and the bound of
Lemma~\ref{lem:k/m} suffices for our purposes.

\begin{lemma}\label{lem:m3}
Assume that a line $\ell$ meets every element of a set $M$ of $m$
lines in three-space such that at most $m/5$ lines of $M$ may be
coplanar or concurrent. Then $\ell$ is incident to at least
$\Omega(m^3)$ distinct reguli spanned by $M$.
\end{lemma}

\begin{proof}
  It is easy to see that there are at least $m^3/50$ (unordered)
  triples of pairwise skew lines in $M$. First notice that there are
  at least $m(3m/5)(m/5)=3m^3/25$ ordered triples
  $(\ell_1,\ell_2,\ell_3)\in M^3$ of pairwise skew lines in $M$:
  Choose any line $\ell_1\in M$; then choose any line $\ell_2\in M$
  that is not incident to the point $\ell\cap \ell_1$ and does not lie
  in the plane $\pi(\ell,\ell_1)$ spanned by $\ell$ and $\ell_1$
  (there are at least $3m/5$ such lines); finally choose any line
  $\ell_3\in M$ that is not incident to $\ell\cap \ell_1$ or $\ell\cap
  \ell_2$ and does not lie in the planes $\pi(\ell,\ell_1)$ or
  $\pi(\ell,\ell_2)$.

  These $m^3/50$ triples of pairwise skew lines do not necessarily span distinct
  reguli, but we show that only few reguli can be incident to too many lines
  of $M$. For every $t\in \N$, let $R_t$ denote the set of reguli spanned
  by $M$ and incident to at least $t$ lines of $M$. There are at
  least $t|R_t|$ incidences on these reguli, and by Lemma~\ref{lem:k/m},
  this number is bounded by $O(m^{3/5}|R_t|^{4/5}+m+|R_t|)$.
  It follows that $|R_t|\leq O(m^3/t^5+m/t)$. Set $t$ to be a large constant
  such that $|R_t|\leq m^3/100$. It follows that at least $m^3/100$ triples
  of $M$ span reguli, each incident to less than $t$ lines. Hence $M$
  spans at least $m^3/(100 {t\choose 3})=\Omega(m^3)$ distinct reguli.
\qed
\end{proof}

\section{Proof of the main theorem\label{sec:proof}}

We are given a set $L$ of $n^2$ lines and a set $S_0$ of points in
$\R^3$ such that every line in $L$ is incident to $n$ points and no
$n$ lines are coplanar. There are $n^3$ point-line incidences, where
each incidence is a pair $(p,\ell)\in S_0\times L$ with $p\in \ell$.
Let $S$ be the set of points in $S_0$ incident to at most $n$ lines
of $L$ and set $s=|S|$. By the Szemer\'edi-Trotter theorem, at most
$O(|L|^2/n^3+|L|/n)=O(n)$ lines of $S_0$ are incident to $n$ or more
lines, and these points are involved in at most $O(n^2\log n)$
incidences.
If $n$ is sufficiently large, then the remaining $s$ points in $S$
and $n^2$ lines in $L$ still have at least $n^3/2$ incidences,
furthermore, no $n$ lines of $L$ are coplanar or meet at a point of
$S$. Let $d=n^3/s$ denote the average number of lines incident to a
point of $S$.

Project the lines of $L$ and the points of $S$ into a generic plane,
and consider the dual arrangement. We obtain a set $L^*$ of $n^2$
points and a set $S^*$ of $s$ lines in the plane such that every point
in $L^*$ is incident to at most $n$ lines, every line in $S^*$ is
incident to at most $n$ points, and there are at least $n^3/2$
point-line incidences.  Choose a parameter $r=cd$ with a sufficiently
small constant $c>0$ to be specified by two upper bounds below.
Consider a {\em $(1/r)$-cutting}~\cite{cf-dvrsu-90} for $S^*$, which
is a partition of the plane into $O(r^2)$ triangles such that the
interior of each triangle intersects at most $s/r$ lines of $S^*$. By
splitting some triangles, if necessary, we obtain a partition of the
plane into a set $\Xi$ of $O(r^2)$ triangles, each containing at most
$n^2/r^2$ points of $L^*$. For every triangle $\sigma\in \Xi$, let
$L^*_\sigma$ denote the set of points of $L^*$ in $\sigma$, and let
$S^*_\sigma$ be the set of lines of $S^*$ intersecting the interior of
$\sigma$.

The number of point-line incidences involving points on the boundary
of some triangles and lines intersecting the interior of an adjacent
triangle is bounded by $I_{\rm boundary}=O(r^2)\cdot
(s/r)=O(sr)=O(csd)=O(cn^3)$. Let $c>0$ be so small that $I_{\rm
  boundary}< n^3/2$. Hence the sum of incidences in each triangle is
at least $n^3/4$, that is,
$$\sum_{\sigma\in \Xi} I(L^*_\sigma,S^*_\sigma) \geq \frac{n^3}{4}.$$
There is a triangle $\sigma$ such that $I(L^*_\sigma,S^*_\sigma)\geq
(n^3/4)/|\Xi|=\Omega(n^3/(c^2d^2))$.  Let $L_\sigma\subset L$ and
$S_\sigma\subset S$ denote, respectively, the lines and points
corresponding to planar duals of $L^*_\sigma$ and $S^*_\sigma$.  We
have a set $L_\sigma$ of at most $O(n^2/(c^2d^2))$ lines and a set
$S_\sigma$ of at most $s/(cd)$ points in three-space that have at
least $\Omega(n^3/(c^2d^2))$ incidences.

We next give a lower bound on the number of line pairs
$$G=\{(\ell_1,\ell_2)\in L_\sigma : \ell_1\cap
\ell_2\in S_\sigma\}$$
that meet at a point of $S_\sigma$. Denoting by $d_\sigma(p)$ the
number of lines of $L_\sigma$ incident to a point $p\in S_\sigma$, we
have
$$|G| = \sum_{p\in S_\sigma} {d_p\choose 2} \geq |S_\sigma|\cdot
{(\sum_{p\in S_\sigma} d_p)/|S_\sigma|\choose 2}.$$
We can estimate the average degree by
$$\frac{\sum_{p\in S_\sigma} d_p}{|S_\sigma|} =
\frac{I(S_\sigma,L_\sigma)}{|S_\sigma|} \geq \Omega
\left(\frac{n^3}{cds}\right) =\Omega \left(\frac{1}{c}\right).$$
Hence, at least $|G| = \Omega(|S_\sigma|/c^2)=\Omega(s/(c^3d))$ line
pairs meet at points of $S_\sigma$. Discard all lines $\ell\in
L_\sigma$ that meet less than $5n$ other lines of $L_\sigma$ at points
of $S_\sigma$.  We have discarded at most
$5n|L_\sigma|=O(n^3/(c^2d^2))= O(s/(c^2d))$ line pairs. Set the
constant $c>0$ so small that we discard at most $|G|/2$ line pairs of
$G$. In the remainder of the proof, $c$ is fixed and hidden in the
asymptotic notation.

We have a set $L_\sigma'$ of $O(n^2/d^2) = O(n^2/(n^3/s)^2) =
O(s^2/n^4)$ lines such that each line meets at least $5n$ other
lines of $L_\sigma'$ and the total number of meeting pairs of lines
is at least $\Omega(s/d)=\Omega(s^2/n^3)$. Recall that no $n$ lines
of $L$ are coplanar or meet at a point of $S$. Let $R_\sigma'$
denote the set of reguli spanned by lines of $L_\sigma'$.  By
Lemma~\ref{lem:m3}, a line $\ell\in L_\sigma'$ that meets $m_\ell$
other lines of $L_\sigma'$, where $m_\ell\geq 5n$, is incident to at
least $\Omega(m_\ell^3)$ reguli of $R_\sigma'$. The total number of
line-reguli incidences in $L_\sigma'\times R_\sigma'$ is bounded
from below by
\begin{eqnarray}
I(L_\sigma',R_\sigma') &=& \sum_{\ell\in L_\sigma'} \Omega(m_\ell^3)
\geq |L_\sigma'| \cdot \Omega \left(\left(\frac{\sum_{\ell\in L_\sigma'}
      m_\ell}{|L_\sigma'|}\right)^3\right) =
    \Omega \left( \frac{1}{|L_\sigma'|^2} \cdot\left(\sum_{\ell\in L_\sigma'}
      m_\ell\right)^3 \right)\nonumber\\
& = &  \Omega\left(\frac{(s^2/n^3)^3}{(s^2/n^4)^2}\right)  =
       \Omega\left(\frac{s^2}{n}\right).\nonumber
\end{eqnarray}
On the other hand, $|L_\sigma'|=O(s^2/n^4)$ lines can span at most
$${|L_\sigma'|\choose 3}= O\left(\frac{s^6}{n^{12}}\right)$$
reguli. By the K\H{o}vari-S\'os-Tur\'an bound, the number of
line-regulus incidences is bounded by
$$I(L_\sigma',R_\sigma')=\Omega(|L_\sigma'|\cdot
|R_\sigma'|^{2/3}+|R_\sigma'|) = O(|L_\sigma'|^3) =
O\left(\frac{s^6}{n^{12}}\right).$$
Comparing the upper and lower bounds on the number of line-regulus
incidences $I(L_\sigma',R_\sigma')$, we have
$$\Omega \left( \frac{s^2}{n}\right) \leq I(L_\sigma',R_\sigma') \leq
O\left(\frac{s^6}{n^{12}}\right),$$
that is, $s=\Omega(n^{11/4})$, as required.  \qed

\end{document}